\documentclass{article}
\usepackage{amsfonts}
\usepackage{amsmath}

\newtheorem{theorem}{Theorem}

\newtheorem{example}[theorem]{Example}

\newtheorem{remark}[theorem]{Remark}

\input{tcilatex}
\begin{document}

\title{Hypersurface family with a common isoasymptotic curve}
\author{Ergin Bayram, Emin Kasap \\
Ondokuz May\i s University}
\maketitle

\begin{abstract}
In the present paper, we handle the problem of finding a hypersurface family
from a given asymptotic curve in $%
\mathbb{R}
^{4}.$ Using the Frenet frame of the given asymptotic curve, we express the
hypersurface as a linear combination of this frame and analyze the necessary
and sufficient conditions for that curve to be asymptotic. We illustrate
this method by presenting some examples.
\end{abstract}



\section{Introduction}

Asymptotic curves are encountered in differential geometry frequently. A
surface curve is called asymptotic if its tangent vectors always point in an
asymptotic direction, that is, the direction in which the normal curvature
is zero. In an asymptotic direction, the surface is not bending away from
its tangent plane.

Asymptotic curves on a surface can be seen in many differential geometry
books $\cite{Oneill}\ -\cite{Struik}.\ $Rastogi \cite{Rastogi} obtained the
differential equation of hyper-asymptotic curves by a new method and showed
some properties of these curves. Aminov \cite{Aminov} established more
general expressions for the curvature of asymptotic curves of submanifolds
in the Riemannian space. Romero - Fuster et.al. \cite{Romero} studied
asymptotic curves on generally immersed surfaces in $%
\mathbb{R}
^{5}.$ Both the general and rational developable surface pencils through an
arbitrary parametric curve as its common asymptotic curve were analyzed by
Liu and Wang \cite{Liu}. Bayram et. al. \cite{Bayram} tackled the problem of
finding a surface pencil from a given asymptotic curve.

However, while differential geometry of a parametric surface in $%
\mathbb{R}
^{3}$ can be found in textbooks such as in Struik \cite{Struik}, Willmore 
\cite{Willmore}, Stoker \cite{Stoker}, do Carmo \cite{Carmo}, differential
geometry of a parametric surface in $%
\mathbb{R}
^{n}$ can be found in textbooks such as in the contemporary literature on
Geometric Modeling \cite{Farin}, \cite{Hoschek}. Also, there is little
literature on differential geometry of parametric surface family in $%
\mathbb{R}
^{3}$ \cite{Alessio}, \cite{Duldul}, \cite{Kasap}, \cite{Wang}, \cite{Bayram}%
, but not in $%
\mathbb{R}
^{4}$. Besides, there is an ascending interest on fourth dimension \cite%
{Abdel}, \cite{Alessio}, \cite{Duldul}.

Furthermore, various visualization techniques about objects in Euclidean
n-space $\left( n\geq 4\right) $ are presented \cite{Banchoff}, \cite%
{Banchoff2}, \cite{Hanson}. The fundamental step to visualize a 4D object is
projecting first into the 3-space and then into the plane. In many real
world applications, the problem of visualizing three-dimensional data,
commonly referred to as scalar fields arouses. The graph of a function $%
\mathbf{f}\left( x,y,z\right) :U\subset 
\mathbb{R}
^{3}\rightarrow 
\mathbb{R}
$, where $U$ is open, is a special type of parametric hypersurface with the
parametrization $\left( x,y,z,\mathbf{f}\left( x,y,z\right) \right) $ in
4-space. There exists a method for rendering such a 3-surface based on known
methods for visualizing functions of two variables \cite{Hamann}.

In this paper, we consider the four dimensional analogue problem of
constructing a parametric representation of a surface family from a given
asymptotic as in Bayram et al. \cite{Bayram}, who derived the necessary and
sufficient conditions on the marching-scale functions for which the curve C
is an asymptotic curve on a given surface. We express the hypersurface
pencil parametrically with the help of the Frenet frame $\left\{ \mathbf{T},%
\mathbf{N},\mathbf{B}_{1},\mathbf{B}_{2}\right\} $ of the given curve. We
find the necessary and sufficient constraints on the marching-scale
functions, namely, coefficients of Frenet vectors, so that both the
asymptotic and parametric requirements are met.

\section{Preliminaries}

Let us first introduce some notations and definitions. Bold letters such as $%
\mathbf{a}$, $\mathbf{R}$ will be used for vectors and vector functions. We
assume that they are smooth enough so that all the (partial) derivatives
given in the paper are meaningful. Let $\mathbf{\alpha :I\subset 
\mathbb{R}
\rightarrow 
\mathbb{R}
}^{4}$ be an arc-length curve. If $\left\{ \mathbf{T},\mathbf{N},\mathbf{B}%
_{1},\mathbf{B}_{2}\right\} $ is the moving Frenet frame along $\mathbf{%
\alpha }$, then the Frenet formulas are given by%
\begin{equation}
\left\{ 
\begin{array}{c}
\mathbf{T}^{\prime }=\kappa _{1}\mathbf{N,} \\ 
\mathbf{N}^{\prime }=-\kappa _{1}\mathbf{T+}\kappa _{2}\mathbf{B}_{1}\mathbf{%
,} \\ 
\mathbf{B}_{1}^{\prime }=-\kappa _{2}\mathbf{N+}\kappa _{3}\mathbf{B}_{2},
\\ 
\mathbf{B}_{2}^{\prime }=-\kappa _{3}\mathbf{B}_{1},%
\end{array}%
\right.  \label{1}
\end{equation}%
where $\mathbf{T},\mathbf{N},\mathbf{B}_{1}$ and $\mathbf{B}_{2}\ $denote
the tangent, principal normal, first binormal and second binormal vector
fields, respectively, $\kappa _{i}\left( i=1,2,3\right) \ $the i-th
curvature functions of the curve $\mathbf{\alpha \ \cite{Hanson}}$.

From elementary differential geometry we have%
\begin{equation}
\left\{ 
\begin{array}{c}
\mathbf{\alpha }^{\prime }\left( s\right) =\mathbf{T}\left( s\right) , \\ 
\mathbf{\alpha }^{\prime \prime }\left( s\right) =\kappa _{1}\left( s\right) 
\mathbf{N}\left( s\right) , \\ 
\kappa _{1}\left( s\right) =\left\Vert \alpha ^{\prime \prime }\left(
s\right) \right\Vert .%
\end{array}%
\right.  \label{2}
\end{equation}

Using Frenet formulas one can obtain the followings%
\begin{equation}
\left\{ 
\begin{array}{c}
\mathbf{\alpha }^{\prime \prime \prime }\left( s\right) =-\kappa _{1}^{2}%
\mathbf{T}\left( s\right) +\kappa _{1}^{\prime }\mathbf{N}\left( s\right)
+\kappa _{1}\kappa _{2}\mathbf{B}_{1}\left( s\right) , \\ 
\mathbf{\alpha }^{\left( iv\right) }\left( s\right) =-3\kappa _{1}\kappa
_{1}^{\prime }\mathbf{T}\left( s\right) +\left( -\kappa _{1}^{3}+\kappa
_{1}^{\prime \prime }-\kappa _{1}\kappa _{2}^{2}\right) \mathbf{N}\left(
s\right) \\ 
+\left( 2\kappa _{1}^{\prime }\kappa _{2}+\kappa _{1}\kappa _{2}^{\prime
}\right) \mathbf{B}_{1}\left( s\right) +\kappa _{1}\kappa _{2}\kappa _{3}%
\mathbf{B}_{2}\left( s\right) .%
\end{array}%
\right.  \label{3}
\end{equation}

The unit vectors $\mathbf{B}_{2}\ $and $\mathbf{B}_{1}\ $are given by%
\begin{equation}
\left\{ 
\begin{array}{c}
\mathbf{B}_{2}\left( s\right) =\frac{\mathbf{\alpha }^{\prime }\left(
s\right) \otimes \mathbf{\alpha }^{\prime \prime }\left( s\right) \otimes 
\mathbf{\alpha }^{\prime \prime \prime }\left( s\right) }{\left\Vert \mathbf{%
\alpha }^{\prime }\left( s\right) \otimes \mathbf{\alpha }^{\prime \prime
}\left( s\right) \otimes \mathbf{\alpha }^{\prime \prime \prime }\left(
s\right) \right\Vert }, \\ 
\mathbf{B}_{1}\left( s\right) =\mathbf{B}_{2}\left( s\right) \otimes \mathbf{%
T}\left( s\right) \otimes \mathbf{N}\left( s\right) ,%
\end{array}%
\right.  \label{4}
\end{equation}%
where $\otimes $ is the vector product of vectors in $%
\mathbb{R}
^{4}$.

Since the vectors $\mathbf{T},\ \mathbf{N},\ \mathbf{B}_{1},\ \mathbf{B}_{2}$
are orthonormal, the second curvature $\kappa _{2}$ and the third curvature $%
\kappa _{3}$ can be obtained from (3) as%
\begin{equation}
\left\{ 
\begin{array}{c}
\kappa _{2}\left( s\right) =\frac{\mathbf{B}_{1}\left( s\right) \bullet 
\mathbf{\alpha }^{\prime \prime \prime }\left( s\right) }{\kappa _{1}\left(
s\right) }, \\ 
\kappa _{3}\left( s\right) =\frac{\mathbf{B}_{2}\left( s\right) \bullet 
\mathbf{\alpha }^{\left( iv\right) }\left( s\right) }{\kappa _{1}\left(
s\right) \kappa _{2}\left( s\right) },%
\end{array}%
\right.  \label{5}
\end{equation}%
where `$\bullet $' denotes the standard inner product.

Let $\left\{ \mathbf{e}_{1},\mathbf{e}_{2},\mathbf{e}_{3},\mathbf{e}%
_{4}\right\} $ be the standard basis for four-dimensional Euclidean space $%
\mathbb{R}
^{4}$. The vector product of the vectors $\mathbf{u=}\sum_{i=1}^{4}u_{i}%
\mathbf{e}_{i},\ \mathbf{v=}\sum_{i=1}^{4}v_{i}\mathbf{e}_{i},$ $\mathbf{w=}%
\sum_{i=1}^{4}w_{i}\mathbf{e}_{i}\ $is defined by%
\begin{equation*}
\mathbf{u\otimes v\otimes w=}\left\vert 
\begin{array}{c}
\mathbf{e}_{1}\ \ \mathbf{e}_{2}\ \ \mathbf{e}_{3}\ \ \mathbf{e}_{4} \\ 
u_{1}\ \ u_{2}\ \ u_{3}\ \ u_{4} \\ 
v_{1}\ \ v_{2}\ \ v_{3}\ \ v_{4} \\ 
w_{1}\ w_{2}\ w_{3}\ w_{4}%
\end{array}%
\right\vert
\end{equation*}

\cite{Hollasch}, \cite{Williams}.

If $\mathbf{u,v}$ and \textbf{w} are linearly independent then $\mathbf{%
u\otimes v\otimes w}$ is orthogonal to each of these vectors.

\section{Hypersurface Family with a Common Isoasymptotic}

A curve $\mathbf{r}\left( s\right) $ on a hypersurface $\mathbf{P=P}\left(
s,t,q\right) \subset 
\mathbb{R}
^{4}$ is called an isoparametric curve if it is a parameter curve, that is,
there exists a pair of parameters $t_{0}$ and $q_{0}$ such that $\mathbf{r}%
\left( s\right) =\mathbf{P}\left( s,t_{0},q_{0}\right) $. Given a parametric
curve $\mathbf{r}\left( s\right) $, it is called an \textit{isoasymptotic}
of a hypersurface $\mathbf{P}$ if it is both a asymptotic and an
isoparametric curve on $\mathbf{P}$.

Let $C:\mathbf{r}=\mathbf{r}\left( s\right) ,\ L_{1}\leq s\leq L_{2}$, be a $%
C^{3}$ curve, where $s$ is the arc-length. To have a well-defined principal
normal, assume that $\mathbf{r}^{\prime \prime }\left( s\right) \neq 0,\
L_{1}\leq s\leq L_{2}$.

Let $\mathbf{T}\left( s\right) ,\mathbf{N}\left( s\right) ,\mathbf{B}%
_{1}\left( s\right) ,\mathbf{B}_{2}\left( s\right) $ be the tangent,
principal normal, first binormal, second binormal, respectively; and let $%
\kappa _{1}\left( s\right) ,\kappa _{2}\left( s\right) $ and $\kappa
_{3}\left( s\right) $ be the first, the second and the third curvature,
respectively. Since $\left\{ \mathbf{T}\left( s\right) ,\mathbf{N}\left(
s\right) ,\mathbf{B}_{1}\left( s\right) ,\mathbf{B}_{2}\left( s\right)
\right\} $ is an orthogonal coordinate frame on $\mathbf{r}\left( s\right) $
the parametric hypersurface $\mathbf{P}\left( s,t,q\right) :\left[
L_{1},L_{2}\right] \times \left[ T_{1},T_{2}\right] \times \left[ Q_{1},Q_{2}%
\right] \rightarrow 
\mathbb{R}
^{4}$ passing through $\mathbf{r}\left( s\right) $ can be defined as follows:%
\begin{equation}
\mathbf{P}\left( s,t,q\right) =\mathbf{r}\left( s\right) +\left( \mathbf{u}%
\left( s,t,q\right) ,\mathbf{v}\left( s,t,q\right) ,\mathbf{w}\left(
s,t,q\right) ,\mathbf{x}\left( s,t,q\right) \right) \left( 
\begin{array}{c}
\mathbf{T}\left( s\right)  \\ 
\mathbf{N}\left( s\right)  \\ 
\mathbf{B}_{1}\left( s\right)  \\ 
\mathbf{B}_{2}\left( s\right) 
\end{array}%
\right) ,  \label{6}
\end{equation}%
\begin{equation*}
L_{1}\leq s\leq L_{2},\ T_{1}\leq s\leq T_{2},\ Q_{1}\leq s\leq Q_{2},
\end{equation*}%
where $\mathbf{u}\left( s,t,q\right) ,\mathbf{v}\left( s,t,q\right) ,\mathbf{%
w}\left( s,t,q\right) $ and $\mathbf{x}\left( s,t,q\right) $ are all $C^{4}$
functions. These functions are called the \textit{marching scale functions}.

We try to find out the necessary and sufficient conditions for which a
hypersurface $\mathbf{P=P}\left( s,t,q\right) $ has the curve $C$ as an
isoasymptotic.

First, to satisfy the isoparametricity condition there should exist $%
t_{0}\in \left[ T_{1},T_{2}\right] $ and $q_{0}\in \left[ Q_{1},Q_{2}\right] 
$ such that $\mathbf{P}\left( s,t_{0},q_{0}\right) =\mathbf{r}\left(
s\right) $, $L_{1}\leq s\leq L_{2}$ , that is,%
\begin{equation}
\left\{ 
\begin{array}{c}
\mathbf{u}\left( s,t_{0},q_{0}\right) =\mathbf{v}\left( s,t_{0},q_{0}\right)
=\mathbf{w}\left( s,t_{0},q_{0}\right) =\mathbf{x}\left(
s,t_{0},q_{0}\right) \equiv 0, \\ 
t_{0}\in \left[ T_{1},T_{2}\right] ,\ q_{0}\in \left[ Q_{1},Q_{2}\right] ,\
L_{1}\leq s\leq L_{2}.%
\end{array}%
\right.   \label{7}
\end{equation}

Secondly, the curve $C$ is an asymptotic curve on the hypersurface $\mathbf{P%
}\left( s,t,q\right) $ if and only if the normal curvature $\kappa
_{n}=S\left( T\right) \bullet T=0$ along the curve, where $S$ is the shape
operator and $T$ is the tangent vector to the curve. The normal $\overset{%
\wedge }{\mathbf{n}}\left( s,t_{0},q_{0}\right) $ of the hypersurface can be
obtained by calculating the vector product of the partial derivatives and
using the Frenet formula as follows

$\frac{\partial \mathbf{P}\left( s,t,q\right) }{\partial s}=\left( 1+\frac{%
\partial \mathbf{u}\left( s,t,q\right) }{\partial s}-\mathbf{v}\left(
s,t,q\right) \kappa _{1}\left( s\right) \right) \mathbf{T}\left( s\right) $

\ \ \ \ \ \ \ \ \ \ $+\left( \mathbf{u}\left( s,t,q\right) \kappa _{1}\left(
s\right) +\frac{\partial \mathbf{v}\left( s,t,q\right) }{\partial s}-\mathbf{%
w}\left( s,t,q\right) \kappa _{2}\left( s\right) \right) \mathbf{N}\left(
s\right) $

\ \ \ \ \ \ \ \ \ \ $+\left( \mathbf{v}\left( s,t,q\right) \kappa _{2}\left(
s\right) +\frac{\partial \mathbf{w}\left( s,t,q\right) }{\partial s}-\mathbf{%
x}\left( s,t,q\right) \kappa _{3}\left( s\right) \right) \mathbf{B}%
_{1}\left( s\right) $

\ \ \ \ \ \ \ \ \ \ $+\left( \mathbf{w}\left( s,t,q\right) \kappa _{3}\left(
s\right) +\frac{\partial \mathbf{x}\left( s,t,q\right) }{\partial s}\right) 
\mathbf{B}_{2}\left( s\right) ,$

\bigskip

$\frac{\partial \mathbf{P}\left( s,t,q\right) }{\partial t}=\frac{\partial 
\mathbf{u}\left( s,t,q\right) }{\partial t}\mathbf{T}\left( s\right) +\frac{%
\partial \mathbf{v}\left( s,t,q\right) }{\partial t}\mathbf{N}\left(
s\right) +\frac{\partial \mathbf{w}\left( s,t,q\right) }{\partial t}\mathbf{B%
}_{1}\left( s\right) +\frac{\partial \mathbf{x}\left( s,t,q\right) }{%
\partial t}\mathbf{B}_{2}\left( s\right) ,$

\bigskip

and

\bigskip

$\frac{\partial \mathbf{P}\left( s,t,q\right) }{\partial q}=\frac{\partial 
\mathbf{u}\left( s,t,q\right) }{\partial q}\mathbf{T}\left( s\right) +\frac{%
\partial \mathbf{v}\left( s,t,q\right) }{\partial q}\mathbf{N}\left(
s\right) +\frac{\partial \mathbf{w}\left( s,t,q\right) }{\partial q}\mathbf{B%
}_{1}\left( s\right) +\frac{\partial \mathbf{x}\left( s,t,q\right) }{%
\partial q}\mathbf{B}_{2}\left( s\right) .$

\begin{remark}
Because,
\end{remark}

\begin{equation*}
\left\{ 
\begin{array}{c}
\mathbf{u}\left( s,t_{0},q_{0}\right) =\mathbf{v}\left( s,t_{0},q_{0}\right)
=\mathbf{w}\left( s,t_{0},q_{0}\right) =\mathbf{x}\left(
s,t_{0},q_{0}\right) \equiv 0, \\ 
t_{0}\in \left[ T_{1},T_{2}\right] ,\ q_{0}\in \left[ Q_{1},Q_{2}\right] ,\
L_{1}\leq s\leq L_{2}.%
\end{array}%
\right. 
\end{equation*}

\textit{along the curve }$\mathit{C}$\textit{, by the definition of partial
differentiation we have}%
\begin{equation*}
\left\{ 
\begin{array}{c}
\frac{\partial \mathbf{u}\left( s,t_{0},q_{0}\right) }{\partial s}=\frac{%
\partial \mathbf{v}\left( s,t_{0},q_{0}\right) }{\partial s}=\frac{\partial 
\mathbf{w}\left( s,t_{0},q_{0}\right) }{\partial s}=\frac{\partial \mathbf{x}%
\left( s,t_{0},q_{0}\right) }{\partial s}\equiv 0, \\ 
t_{0}\in \left[ T_{1},T_{2}\right] ,\ q_{0}\in \left[ Q_{1},Q_{2}\right] ,\
L_{1}\leq s\leq L_{2}.%
\end{array}%
\right. 
\end{equation*}

\qquad Using $\left( \ref{7}\right) $ we have

$\overset{\wedge }{\mathbf{n}}\left( s,t_{0},q_{0}\right) =\frac{\partial 
\mathbf{P}\left( s,t_{0},q_{0}\right) }{\partial s}\otimes \frac{\partial 
\mathbf{P}\left( s,t_{0},q_{0}\right) }{\partial t}\otimes \frac{\partial 
\mathbf{P}\left( s,t_{0},q_{0}\right) }{\partial q}$

\ \ \ \ \ \ \ \ \ \ \ \ \ \ \ $=\phi _{1}\left( s,t_{0},q_{0}\right) \mathbf{%
T}\left( s\right) -\phi _{2}\left( s,t_{0},q_{0}\right) \mathbf{N}\left(
s\right) $

$\ \ \ \ \ \ \ \ \ \ \ \ \ \ \ +\phi _{3}\left( s,t_{0},q_{0}\right) \mathbf{%
B}_{1}\left( s\right) -\phi _{4}\left( s,t_{0},q_{0}\right) \mathbf{B}%
_{2}\left( s\right) ,$

where

$\phi _{1}\left( s,t_{0},q_{0}\right) =\left\vert 
\begin{array}{c}
\frac{\partial \mathbf{v}\left( s,t_{0},q_{0}\right) }{\partial s}\ \frac{%
\partial \mathbf{w}\left( s,t_{0},q_{0}\right) }{\partial s}\ \frac{\partial 
\mathbf{x}\left( s,t_{0},q_{0}\right) }{\partial s} \\ 
\frac{\partial \mathbf{v}\left( s,t_{0},q_{0}\right) }{\partial t}\ \frac{%
\partial \mathbf{w}\left( s,t_{0},q_{0}\right) }{\partial t}\ \frac{\partial 
\mathbf{x}\left( s,t_{0},q_{0}\right) }{\partial t} \\ 
\frac{\partial \mathbf{v}\left( s,t_{0},q_{0}\right) }{\partial q}\ \frac{%
\partial \mathbf{w}\left( s,t_{0},q_{0}\right) }{\partial q}\ \frac{\partial 
\mathbf{x}\left( s,t_{0},q_{0}\right) }{\partial q}%
\end{array}%
\right\vert =0,$

\bigskip

$\phi _{2}\left( s,t_{0},q_{0}\right) =\left\vert 
\begin{array}{c}
1+\frac{\partial \mathbf{u}\left( s,t_{0},q_{0}\right) }{\partial s}\ \frac{%
\partial \mathbf{w}\left( s,t_{0},q_{0}\right) }{\partial s}\ \frac{\partial 
\mathbf{x}\left( s,t_{0},q_{0}\right) }{\partial s} \\ 
\frac{\partial \mathbf{u}\left( s,t_{0},q_{0}\right) }{\partial t}\ \frac{%
\partial \mathbf{w}\left( s,t_{0},q_{0}\right) }{\partial t}\ \frac{\partial 
\mathbf{x}\left( s,t_{0},q_{0}\right) }{\partial t} \\ 
\frac{\partial \mathbf{u}\left( s,t_{0},q_{0}\right) }{\partial q}\ \frac{%
\partial \mathbf{w}\left( s,t_{0},q_{0}\right) }{\partial q}\ \frac{\partial 
\mathbf{x}\left( s,t_{0},q_{0}\right) }{\partial q}%
\end{array}%
\right\vert $

\ \ \ \ \ \ \ \ \ \ \ \ \ \ \ \ $=\left\vert 
\begin{array}{c}
1\ \ \ \ \ \ \ \ \ \ \ \ \ \ 0\ \ \ \ \ \ \ \ \ \ \ \ \ \ 0 \\ 
\frac{\partial \mathbf{u}\left( s,t_{0},q_{0}\right) }{\partial t}\ \frac{%
\partial \mathbf{w}\left( s,t_{0},q_{0}\right) }{\partial t}\ \frac{\partial 
\mathbf{x}\left( s,t_{0},q_{0}\right) }{\partial t} \\ 
\frac{\partial \mathbf{u}\left( s,t_{0},q_{0}\right) }{\partial q}\ \frac{%
\partial \mathbf{w}\left( s,t_{0},q_{0}\right) }{\partial q}\ \frac{\partial 
\mathbf{x}\left( s,t_{0},q_{0}\right) }{\partial q}%
\end{array}%
\right\vert $

\ \ \ \ \ \ \ \ \ \ \ \ \ \ \ \ $=\frac{\partial \mathbf{w}\left(
s,t_{0},q_{0}\right) }{\partial t}\ \frac{\partial \mathbf{x}\left(
s,t_{0},q_{0}\right) }{\partial q}-\frac{\partial \mathbf{w}\left(
s,t_{0},q_{0}\right) }{\partial q}\frac{\partial \mathbf{x}\left(
s,t_{0},q_{0}\right) }{\partial t},$

\bigskip

$\phi _{3}\left( s,t_{0},q_{0}\right) =\left\vert 
\begin{array}{c}
1+\frac{\partial \mathbf{u}\left( s,t_{0},q_{0}\right) }{\partial s}\ \frac{%
\partial \mathbf{v}\left( s,t_{0},q_{0}\right) }{\partial s}\ \frac{\partial 
\mathbf{x}\left( s,t_{0},q_{0}\right) }{\partial s} \\ 
\frac{\partial \mathbf{u}\left( s,t_{0},q_{0}\right) }{\partial t}\ \frac{%
\partial \mathbf{v}\left( s,t_{0},q_{0}\right) }{\partial t}\ \frac{\partial 
\mathbf{x}\left( s,t_{0},q_{0}\right) }{\partial t} \\ 
\frac{\partial \mathbf{u}\left( s,t_{0},q_{0}\right) }{\partial q}\ \frac{%
\partial \mathbf{v}\left( s,t_{0},q_{0}\right) }{\partial q}\ \frac{\partial 
\mathbf{x}\left( s,t_{0},q_{0}\right) }{\partial q}%
\end{array}%
\right\vert $

\ \ \ \ \ \ \ \ \ \ \ \ \ \ \ \ $=\left\vert 
\begin{array}{c}
1\ \ \ \ \ \ \ \ \ \ \ \ \ \ 0\ \ \ \ \ \ \ \ \ \ \ \ \ \ 0 \\ 
\frac{\partial \mathbf{u}\left( s,t_{0},q_{0}\right) }{\partial t}\ \frac{%
\partial \mathbf{v}\left( s,t_{0},q_{0}\right) }{\partial t}\ \frac{\partial 
\mathbf{x}\left( s,t_{0},q_{0}\right) }{\partial t} \\ 
\frac{\partial \mathbf{u}\left( s,t_{0},q_{0}\right) }{\partial q}\ \frac{%
\partial \mathbf{v}\left( s,t_{0},q_{0}\right) }{\partial q}\ \frac{\partial 
\mathbf{x}\left( s,t_{0},q_{0}\right) }{\partial q}%
\end{array}%
\right\vert $

\ \ \ \ \ \ \ \ \ \ \ \ \ \ \ \ $=\frac{\partial \mathbf{v}\left(
s,t_{0},q_{0}\right) }{\partial t}\ \frac{\partial \mathbf{x}\left(
s,t_{0},q_{0}\right) }{\partial q}-\frac{\partial \mathbf{v}\left(
s,t_{0},q_{0}\right) }{\partial q}\frac{\partial \mathbf{x}\left(
s,t_{0},q_{0}\right) }{\partial t},$

\bigskip

\bigskip

$\phi _{4}\left( s,t_{0},q_{0}\right) =\left\vert 
\begin{array}{c}
1+\frac{\partial \mathbf{u}\left( s,t_{0},q_{0}\right) }{\partial s}\ \frac{%
\partial \mathbf{v}\left( s,t_{0},q_{0}\right) }{\partial s}\ \frac{\partial 
\mathbf{w}\left( s,t_{0},q_{0}\right) }{\partial s} \\ 
\frac{\partial \mathbf{u}\left( s,t_{0},q_{0}\right) }{\partial t}\ \frac{%
\partial \mathbf{v}\left( s,t_{0},q_{0}\right) }{\partial t}\ \frac{\partial 
\mathbf{w}\left( s,t_{0},q_{0}\right) }{\partial t} \\ 
\frac{\partial \mathbf{u}\left( s,t_{0},q_{0}\right) }{\partial q}\ \frac{%
\partial \mathbf{v}\left( s,t_{0},q_{0}\right) }{\partial q}\ \frac{\partial 
\mathbf{w}\left( s,t_{0},q_{0}\right) }{\partial q}%
\end{array}%
\right\vert $

\ \ \ \ \ \ \ \ \ \ \ \ \ \ \ \ $=\left\vert 
\begin{array}{c}
1\ \ \ \ \ \ \ \ \ \ \ \ \ \ 0\ \ \ \ \ \ \ \ \ \ \ \ \ \ 0 \\ 
\frac{\partial \mathbf{u}\left( s,t_{0},q_{0}\right) }{\partial t}\ \frac{%
\partial \mathbf{v}\left( s,t_{0},q_{0}\right) }{\partial t}\ \frac{\partial 
\mathbf{w}\left( s,t_{0},q_{0}\right) }{\partial t} \\ 
\frac{\partial \mathbf{u}\left( s,t_{0},q_{0}\right) }{\partial q}\ \frac{%
\partial \mathbf{v}\left( s,t_{0},q_{0}\right) }{\partial q}\ \frac{\partial 
\mathbf{w}\left( s,t_{0},q_{0}\right) }{\partial q}%
\end{array}%
\right\vert $

\ \ \ \ \ \ \ \ \ \ \ \ \ \ \ \ $=\frac{\partial \mathbf{v}\left(
s,t_{0},q_{0}\right) }{\partial t}\ \frac{\partial \mathbf{w}\left(
s,t_{0},q_{0}\right) }{\partial q}-\frac{\partial \mathbf{v}\left(
s,t_{0},q_{0}\right) }{\partial q}\frac{\partial \mathbf{w}\left(
s,t_{0},q_{0}\right) }{\partial t}.$

\bigskip

\bigskip

So,

\begin{eqnarray}
\kappa _{n} &=&S\left( T\right) \bullet T=0\ \ \Leftrightarrow \ \overset{%
\wedge }{\mathbf{n}}\bullet \mathbf{N=0\ \Leftrightarrow }\   \label{8} \\
\phi _{2}\left( s,t_{0},q_{0}\right)  &\equiv &0,\ \phi _{3}^{2}\left(
s,t_{0},q_{0}\right) +\phi _{4}^{2}\left( s,t_{0},q_{0}\right) \neq 0, 
\notag \\
t_{0} &\in &\left[ T_{1},T_{2}\right] ,\ q_{0}\in \left[ Q_{1},Q_{2}\right]
,\ L_{1}\leq s\leq L_{2}.  \notag
\end{eqnarray}

Thus, any hypersurface defined by (\ref{6}) has the curve $C$ as an
isoasymptotic if and only if%
\begin{equation}
\left\{ 
\begin{array}{c}
\mathbf{u}\left( s,t_{0},q_{0}\right) =\mathbf{v}\left( s,t_{0},q_{0}\right)
=\mathbf{w}\left( s,t_{0},q_{0}\right) =\mathbf{x}\left(
s,t_{0},q_{0}\right) \equiv 0, \\ 
\ \phi _{2}\left( s,t_{0},q_{0}\right) \equiv 0,\ \phi _{3}^{2}\left(
s,t_{0},q_{0}\right) +\phi _{4}^{2}\left( s,t_{0},q_{0}\right) \neq 0,%
\end{array}%
\right.   \label{9}
\end{equation}%
\begin{equation*}
t_{0}\in \left[ T_{1},T_{2}\right] ,\ q_{0}\in \left[ Q_{1},Q_{2}\right] ,\
L_{1}\leq s\leq L_{2}.
\end{equation*}

\bigskip is satisfied. We call the set of hypersurfaces defined by (\ref{6})
and satisfying (\ref{9}) an\textit{\ isoasymptotic hypersurface family}.

\section{Examples}

\begin{example}
\bigskip Let $\mathbf{r}\left( s\right) =\left( \frac{1}{2}\cos \left(
s\right) ,\frac{1}{2}\sin \left( s\right) ,\frac{1}{2}s,\frac{\sqrt{2}}{2}%
s\right) ,\ 0\leq s\leq 2\pi ,\ $be a curve parametrized by arc-length. For
this curve,%
\begin{eqnarray*}
\mathbf{T}\left( s\right)  &=&\mathbf{r}^{\prime }\left( s\right) =\left( -%
\frac{1}{2}\sin \left( s\right) ,\frac{1}{2}\cos \left( s\right) ,\frac{1}{2}%
,\frac{\sqrt{2}}{2}\right) , \\
\mathbf{N}\left( s\right)  &=&\left( -\cos \left( s\right) ,-\sin \left(
s\right) ,0,0\right) , \\
\mathbf{B}_{2}\left( s\right)  &=&\frac{\mathbf{r}^{\prime }\left( s\right)
\otimes \mathbf{r}^{\prime \prime }\left( s\right) \otimes \mathbf{r}%
^{\prime \prime \prime }\left( s\right) }{\left\Vert \mathbf{r}^{\prime
}\left( s\right) \otimes \mathbf{r}^{\prime \prime }\left( s\right) \otimes 
\mathbf{r}^{\prime \prime \prime }\left( s\right) \right\Vert }=\left( 0,0,%
\frac{\sqrt{6}}{3},-\frac{\sqrt{3}}{3}\right) , \\
\mathbf{B}_{1}\left( s\right)  &=&\mathbf{B}_{2}\otimes \mathbf{T}\otimes 
\mathbf{N}=\left( -\frac{\sqrt{3}}{2}\sin \left( s\right) ,\frac{\sqrt{3}}{2}%
\cos \left( s\right) ,-\frac{\sqrt{3}}{6},-\frac{\sqrt{6}}{6}\right) .
\end{eqnarray*}%
Let us choose the marching-scale functions as%
\begin{eqnarray*}
\mathbf{u}\left( s,t,q\right)  &=&\left( t-t_{0}\right) \left(
q-q_{0}\right) ,\ \mathbf{v}\left( s,t,q\right) =t-t_{0},\ \mathbf{w}\left(
s,t,q\right) \equiv 0,\ \mathbf{x}\left( s,t,q\right) =q-q_{0},\  \\
t_{0} &\in &\left[ 0,1\right] ,\ q_{0}\in \left[ 0,1\right] ,\ 0\leq s\leq
2\pi .
\end{eqnarray*}%
So, we have the hypersurface%
\begin{eqnarray*}
\mathbf{P}\left( s,t,q\right)  &=&\mathbf{r}\left( s\right) +\mathbf{u}%
\left( s,t,q\right) \mathbf{T}\left( s\right) +\mathbf{v}\left( s,t,q\right) 
\mathbf{N}\left( s\right) +\mathbf{w}\left( s,t,q\right) \mathbf{B}%
_{1}\left( s\right) +\mathbf{x}\left( s,t,q\right) \mathbf{B}_{2}\left(
s\right)  \\
&=&\left( \frac{1}{2}\cos \left( s\right) -\frac{1}{2}\left( t-t_{0}\right)
\left( q-q_{0}\right) \sin \left( s\right) -\left( t-t_{0}\right) \cos
\left( s\right) ,\right.  \\
&&\frac{1}{2}\sin \left( s\right) +\frac{1}{2}\left( t-t_{0}\right) \left(
q-q_{0}\right) \cos \left( s\right) -\left( t-t_{0}\right) \sin \left(
s\right) , \\
&&\frac{1}{2}s+\frac{1}{2}\left( t-t_{0}\right) \left( q-q_{0}\right) +\frac{%
\sqrt{6}}{3}\left( q-q_{0}\right) , \\
&&\left. \frac{\sqrt{2}}{2}s+\frac{\sqrt{2}}{2}\left( t-t_{0}\right) \left(
q-q_{0}\right) -\frac{\sqrt{3}}{3}\left( q-q_{0}\right) \right) ,
\end{eqnarray*}%
$0\leq s\leq 2\pi ,\ 0\leq t\leq 1,\ 0\leq q\leq 1,\ t_{0}\in \left[ 0,1%
\right] ,\ q_{0}\in \left[ 0,1\right] ,\ $is a member of the isoasymptotic
hypersurface family, since it satisfies (\ref{9}). \newline
By changing the parameters $t_{0}\ $and $q_{0}\ $we can adjust the position
of the curve $\mathbf{r}\left( s\right) $ on the hypersurface. Let us choose 
$t_{0}=\frac{1}{2}$ and $q_{0}=0$. Now the curve $\mathbf{r}\left( s\right) $
is again an isoasymptotic on the hypersurface $\mathbf{P}\left( s,t,q\right) 
$ and the equation of the hypersurface is%
\begin{eqnarray*}
\mathbf{P}\left( s,t,q\right)  &=&\left( \cos s-t\cos s+\frac{1}{4}q\sin s-%
\frac{1}{2}qt\sin s,\right.  \\
&&\sin s-\frac{1}{4}q\cos s-t\sin s+\frac{1}{2}qt\cos s, \\
&&\frac{1}{2}s-\frac{1}{4}q+\frac{1}{3}\sqrt{6}q+\frac{1}{2}qt, \\
&&\left. \frac{1}{2}\sqrt{2}s-\frac{1}{3}\sqrt{3}q-\frac{1}{4}\sqrt{2}q+%
\frac{1}{2}\sqrt{2}qt\right) .
\end{eqnarray*}%
The projection of a hypersurface into 3-space generally yields a
three-dimensional volume. If we fix each of the three parameters, one at a
time, we obtain three distinct families of 2-spaces in 4-space. The
projections of these 2-surfaces into 3-space are surfaces in 3-space. Thus,
they can be displayed by 3D rendering methods. \qquad \newline
So, if we (parallel) project the hypersurface $\mathbf{P}\left( s,t,q\right) 
$ into the $\mathbf{w}=0$ subspace and fix $q=0$ we obtain the surface%
\begin{equation*}
\mathbf{P}_{\mathbf{w}}\left( s,t,0\right) =\left( \cos s-t\cos s,\sin
s-t\sin s,\frac{1}{2}s\right) ,
\end{equation*}%
$0\leq s\leq 2\pi ,0\leq t\leq 1\ $in 3-space illustrated in Fig. 1.\FRAME{%
dtbpFUX}{4.1563in}{2.7709in}{0pt}{\Qcb{Fig. 1. Projection of a member of the
hypersurface family and its isoasymptotic.}}{}{Plot}{\special{language
"Scientific Word";type "MAPLEPLOT";width 4.1563in;height 2.7709in;depth
0pt;display "USEDEF";plot_snapshots TRUE;mustRecompute FALSE;lastEngine
"MuPAD";xmin "0";xmax "6.28";ymin "0";ymax "1";xviewmin "-0.999999";xviewmax
"1";yviewmin "-0.999997";yviewmax "1";zviewmin "-3.14E-10";zviewmax
"3.14";phi 92;theta -176;cameraDistance "16.0292";cameraOrientation
"[0,0,0.212314]";cameraOrientationFixed TRUE;plottype 5;axesFont "Times New
Roman,12,0000000000,useDefault,normal";num-x-gridlines 25;num-y-gridlines
25;plotstyle "wireframe";axesstyle "none";axestips FALSE;plotshading
"XYZ";lighting 0;xis \TEXUX{x};yis \TEXUX{y};var1name \TEXUX{$x$};var2name
\TEXUX{$y$};function \TEXUX{$\left( \cos s-t\cos s,\sin s-t\sin
s,\frac{1}{2}s\right) $};linestyle 1;pointstyle "point";linethickness
1;lineAttributes "Solid";var1range "0,6.28";var2range "0,1";surfaceColor
"[linear:XYZ:RGB:0x0000ff00:0x0000ff00]";surfaceStyle "Wire
Frame";num-x-gridlines 25;num-y-gridlines 25;surfaceMesh
"Mesh";rangeset"XY";function \TEXUX{$\left( \frac{1}{2}\cos \left( s\right)
,\frac{1}{2}\sin \left( s\right) ,\frac{1}{2}s\right) $};linestyle
1;pointstyle "point";linethickness 3;lineAttributes "Solid";curveStyle
"Line";var1range "0,6.28";var2range "-5,5";surfaceColor
"[linear:XYZ:RGB:0000000000:0000000000]";num-x-gridlines 25;num-y-gridlines
25;rangeset"X";VCamFile 'NAYE4J1J.xvz';valid_file "T";tempfilename
'NAYE4J0C.wmf';tempfile-properties "XPR";}}
\end{example}

\begin{example}
Given the curve $\mathbf{r}\left( s\right) =\left( \frac{1}{2}\sin \left(
s\right) ,\frac{1}{2}\cos \left( s\right) ,0,\frac{\sqrt{3}}{2}s\right) ,\
0\leq s\leq 3,$ it is easy to show that%
\begin{eqnarray*}
\mathbf{T}\left( s\right)  &=&\mathbf{r}^{\prime }\left( s\right) =\left( 
\frac{1}{2}\cos \left( s\right) ,-\frac{1}{2}\sin \left( s\right) ,0,\frac{%
\sqrt{3}}{2}\right) , \\
\mathbf{N}\left( s\right)  &=&\left( -\sin \left( s\right) ,-\cos \left(
s\right) ,0,0\right) , \\
\mathbf{B}_{2}\left( s\right)  &=&\frac{\mathbf{r}^{\prime }\left( s\right)
\otimes \mathbf{r}^{\prime \prime }\left( s\right) \otimes \mathbf{r}%
^{\prime \prime \prime }\left( s\right) }{\left\Vert \mathbf{r}^{\prime
}\left( s\right) \otimes \mathbf{r}^{\prime \prime }\left( s\right) \otimes 
\mathbf{r}^{\prime \prime \prime }\left( s\right) \right\Vert }=\left(
0,0,-1,0\right) , \\
\mathbf{B}_{1}\left( s\right)  &=&\mathbf{B}_{2}\otimes \mathbf{T}\otimes 
\mathbf{N}=\left( \frac{\sqrt{3}}{2}\cos \left( s\right) ,-\frac{\sqrt{3}}{2}%
\sin \left( s\right) ,0,-\frac{1}{2}\right) .
\end{eqnarray*}%
Let us choose the marching-scale functions as%
\begin{eqnarray*}
\mathbf{u}\left( s,t,q\right)  &=&\left( t-t_{0}\right) , \\
\mathbf{v}\left( s,t,q\right)  &=&\left( s+t+1\right) \left( q-q_{0}\right) ,
\\
\mathbf{w}\left( s,t,q\right)  &\equiv &0, \\
\mathbf{x}\left( s,t,q\right)  &=&\left( s+1\right) \left( t-t_{0}\right) .
\end{eqnarray*}%
From (\ref{9}), the hypersurface%
\begin{eqnarray*}
\mathbf{P}\left( s,t,q\right)  &=&\mathbf{r}\left( s\right) +\mathbf{u}%
\left( s,t,q\right) \mathbf{T}\left( s\right) +\mathbf{v}\left( s,t,q\right) 
\mathbf{N}\left( s\right) +\mathbf{w}\left( s,t,q\right) \mathbf{B}%
_{1}\left( s\right) +\mathbf{x}\left( s,t,q\right) \mathbf{B}_{2}\left(
s\right)  \\
&=&\left( \frac{1}{2}\sin \left( s\right) -\left( s+t+1\right) \left(
q-q_{0}\right) \sin \left( s\right) +\frac{1}{2}\left( t-t_{0}\right) \cos
\left( s\right) ,\right.  \\
&&\frac{1}{2}\cos \left( s\right) -\left( s+t+1\right) \left( q-q_{0}\right)
\cos \left( s\right) -\frac{1}{2}\left( t-t_{0}\right) \sin \left( s\right) ,
\\
&&-\left( s+1\right) \left( t-t_{0}\right) , \\
&&\left. \frac{\sqrt{3}}{2}s+\frac{\sqrt{3}}{2}\left( t-t_{0}\right) \right)
,
\end{eqnarray*}%
$0\leq s\leq 3,\ 0\leq t\leq 1,\ 0\leq q\leq 1,\ $is a member of the
hypersurface family having the curve $\mathbf{r}\left( s\right) $ as an
isoasymptotic.\newline
\newline
Setting $t_{0}=\frac{1}{2}$ and $q_{0}=0$ yields the hypersurface%
\begin{eqnarray*}
\mathbf{P}\left( s,t,q\right)  &=&\left( \frac{1}{2}\sin \left( s\right)
-\left( s+t+1\right) q\sin \left( s\right) +\frac{1}{2}\left( t-\frac{1}{2}%
\right) \cos \left( s\right) ,\right.  \\
&&\frac{1}{2}\cos \left( s\right) -\left( s+t+1\right) q\cos \left( s\right)
-\frac{1}{2}\left( t-\frac{1}{2}\right) \sin \left( s\right) , \\
&&\left. -\left( s+1\right) \left( t-\frac{1}{2}\right) ,\frac{\sqrt{3}}{2}s+%
\frac{\sqrt{3}}{2}\left( t-\frac{1}{2}\right) \right) ,
\end{eqnarray*}%
By (parallel) projecting the hypersurface $\mathbf{P}\left( s,t,q\right) $
into the subspace $\mathbf{w=0}$ and fixing $q=0$ we get the surface%
\begin{eqnarray*}
\mathbf{P}_{\mathbf{w}}\left( s,t,0\right)  &=&\left( \frac{1}{2}\sin \left(
s\right) +\frac{1}{2}\left( t-\frac{1}{2}\right) \cos \left( s\right)
,\right.  \\
&&\frac{1}{2}\cos \left( s\right) -\frac{1}{2}\left( t-\frac{1}{2}\right)
\sin \left( s\right) , \\
&&\left. -\left( s+1\right) \left( t-\frac{1}{2}\right) \right) ,
\end{eqnarray*}%
where, $0\leq s\leq 3,\ 0\leq t\leq 1$ in 3-space, illustrated in Fig. 2.%
\FRAME{dtbpFUX}{4.1554in}{2.77in}{0pt}{\Qcb{Fig. 2. Projection of a member
of the hypersurface family and its isoasymptotic.}}{}{Plot}{\special%
{language "Scientific Word";type "MAPLEPLOT";width 4.1554in;height
2.77in;depth 0pt;display "USEDEF";plot_snapshots TRUE;mustRecompute
FALSE;lastEngine "MuPAD";xmin "0";xmax "3";ymin "0";ymax "1";xviewmin
"-0.25";xviewmax "0.558928";yviewmin "-0.558234";yviewmax
"0.558648";zviewmin "-2";zviewmax "2";phi 91;theta 3;cameraDistance
"70.8556";cameraOrientation "[0,0,0.0789505]";cameraOrientationFixed
TRUE;plottype 5;axesFont "Times New
Roman,12,0000000000,useDefault,normal";num-x-gridlines 25;num-y-gridlines
25;plotstyle "wireframe";axesstyle "none";axestips FALSE;plotshading
"XYZ";lighting 0;xis \TEXUX{x};yis \TEXUX{y};var1name \TEXUX{$x$};var2name
\TEXUX{$y$};function \TEXUX{$\left( \frac{1}{2}\sin \left( s\right) +\left(
t-\frac{1}{2}\right) \frac{1}{2}\cos \left( s\right) ,\frac{1}{2}\cos \left(
s\right) -\left( t-\frac{1}{2}\right) \frac{1}{2}\sin \left( s\right)
,-\left( s+1\right) \left( t-\frac{1}{2}\right) \right) $};linestyle
1;pointstyle "point";linethickness 1;lineAttributes "Solid";var1range
"0,3";var2range "0,1";surfaceColor
"[linear:XYZ:RGB:0x990000ff:0x990000ff]";surfaceStyle "Wire
Frame";num-x-gridlines 25;num-y-gridlines 25;surfaceMesh
"Mesh";rangeset"XY";function \TEXUX{$\left( \frac{1}{2}\sin \left( s\right)
,\frac{1}{2}\cos \left( s\right) ,0\right) $};linestyle 1;pointstyle
"point";linethickness 3;lineAttributes "Solid";curveStyle "Line";var1range
"0,3";var2range "0,1";surfaceColor
"[linear:XYZ:RGB:0000000000:0000000000]";num-x-gridlines 25;num-y-gridlines
25;rangeset"X";VCamFile 'NAYE4J1L.xvz';valid_file "T";tempfilename
'NAYE4J0D.wmf';tempfile-properties "XPR";}}For the same curve in question
let us choose marching-scale functions as 
\begin{eqnarray*}
\mathbf{u}\left( s,t,q\right)  &\equiv &0, \\
\mathbf{v}\left( s,t,q\right)  &=&\sin \left( s\left( q-q_{0}\right) \right)
, \\
\mathbf{w}\left( s,t,q\right)  &\equiv &0, \\
\mathbf{x}\left( s,t,q\right)  &=&sq^{2}\left( t-t_{0}\right) .
\end{eqnarray*}%
Thus, from (\ref{9}) the curve $\mathbf{r}\left( s\right) $ is an
isoasymptotic on the hypersurface%
\begin{eqnarray*}
\mathbf{P}\left( s,t,q\right)  &=&\mathbf{r}\left( s\right) +\mathbf{u}%
\left( s,t,q\right) \mathbf{T}\left( s\right) +\mathbf{v}\left( s,t,q\right) 
\mathbf{N}\left( s\right) +\mathbf{w}\left( s,t,q\right) \mathbf{B}%
_{1}\left( s\right) +\mathbf{x}\left( s,t,q\right) \mathbf{B}_{2}\left(
s\right)  \\
&=&\left( \frac{1}{2}\sin \left( s\right) -\sin \left( s\right) \sin \left(
s\left( q-q_{0}\right) \right) ,\right.  \\
&&\frac{1}{2}\cos \left( s\right) -\cos \left( s\right) \sin \left( s\left(
q-q_{0}\right) \right) , \\
&&\left. -sq^{2}\left( t-t_{0}\right) ,\frac{\sqrt{3}}{2}s\right) ,
\end{eqnarray*}%
where $0<s\leq \frac{\pi }{2},\ 0\leq t\leq 1,\ 0<q<1.$\newline
\newline
By taking $t_{0}=1$ and $q_{0}=\frac{1}{2}$ we have the following
hypersurface:%
\begin{eqnarray*}
\mathbf{P}\left( s,t,q\right)  &=&\left( \frac{1}{2}\sin \left( s\right)
-\sin \left( s\right) \sin \left( s\left( q-\frac{1}{2}\right) \right)
,\right.  \\
&&\frac{1}{2}\cos \left( s\right) -\cos \left( s\right) \sin \left( s\left(
q-\frac{1}{2}\right) \right) , \\
&&\left. -sq^{2}\left( t-1\right) ,\frac{\sqrt{3}}{2}s\right) .
\end{eqnarray*}%
Hence, if we (parallel) project the hypersurface $\mathbf{P}\left(
s,t,q\right) $ into the $\mathbf{z=0}$ subspace we get the surface%
\begin{eqnarray*}
\mathbf{P}_{z}\left( s,q\right)  &=&\left( \frac{1}{2}\sin \left( s\right)
-\sin \left( s\right) \sin \left( s\left( q-\frac{1}{2}\right) \right)
,\right.  \\
&&\frac{1}{2}\cos \left( s\right) -\cos \left( s\right) \sin \left( s\left(
q-\frac{1}{2}\right) \right) , \\
&&\left. \frac{\sqrt{3}}{2}s\right) ,
\end{eqnarray*}%
where $0<s\leq \frac{\pi }{2},\ 0<q<1,\ $in 3-space shown in Fig. 3.\FRAME{%
dtbpFUX}{4.4996in}{3in}{0pt}{\Qcb{Fig. 3. Projection of a member of the
hypersurface family and its isoasymptotic.}}{}{Plot}{\special{language
"Scientific Word";type "MAPLEPLOT";width 4.4996in;height 3in;depth
0pt;display "USEDEF";plot_snapshots TRUE;mustRecompute FALSE;lastEngine
"MuPAD";xmin "0";xmax "1";ymin "0";ymax "1.55";xviewmin "-0.199673";xviewmax
"1.19946";yviewmin "-0.0280655";yviewmax "0.657381";zviewmin
"-1.34234E-10";zviewmax "1.34234";phi 26;theta -47;cameraDistance
"6.25566";cameraOrientation "[0,0,0.7698]";cameraOrientationFixed
TRUE;plottype 5;axesFont "Times New
Roman,12,0000000000,useDefault,normal";num-x-gridlines 25;num-y-gridlines
25;plotstyle "wireframe";axesstyle "none";axestips FALSE;plotshading
"XYZ";lighting 0;xis \TEXUX{x};yis \TEXUX{y};var1name \TEXUX{$x$};var2name
\TEXUX{$y$};function \TEXUX{$\left( \frac{1}{2}\sin \left( s\right) -\sin
\left( s\right) \sin \left( s\left( q-\frac{1}{2}\right) \right)
,\frac{1}{2}\cos \left( s\right) -\cos \left( s\right) \sin \left( s\left(
q-\frac{1}{2}\right) \right) ,\frac{\sqrt{3}}{2}s\right) $};linestyle
1;pointstyle "point";linethickness 1;lineAttributes "Solid";var1range
"0,1";var2range "0,1.55";surfaceColor
"[linear:XYZ:RGB:0x00ff0000:0x00ff0000]";surfaceStyle "Wire
Frame";num-x-gridlines 25;num-y-gridlines 25;surfaceMesh
"Mesh";rangeset"XY";function \TEXUX{$\left( \frac{1}{2}\sin \left( s\right)
,\frac{1}{2}\cos \left( s\right) ,\frac{\sqrt{3}}{2}s\right) $};linestyle
1;pointstyle "point";linethickness 3;lineAttributes "Solid";curveStyle
"Line";var1range "0,1.55";var2range "0,1";surfaceColor
"[linear:XYZ:RGB:0000000000:0000000000]";num-x-gridlines 25;num-y-gridlines
25;rangeset"X";VCamFile 'NAYE4J1K.xvz';valid_file "T";tempfilename
'NAYE4J0E.wmf';tempfile-properties "XPR";}}
\end{example}

\end{document}